\documentclass[12pt,thmsa]{article}
\usepackage{amssymb}

\input{tcilatex}

\begin{document}

\begin{center}
\textbf{Direct products and elementary equivalence}

\textbf{of polycyclic-by-finite groups}

C. Lasserre and F. Oger\bigskip
\end{center}

\bigskip 

\bigskip 

\bigskip 

\noindent Lasserre Cl\'{e}ment: Institut Fourier - U.M.R. 5582, 100 rue des
Maths, BP 74, 38402 Saint Martin d'H\`{e}res, France;
lasserre.clem@gmail.com.

\noindent Oger Francis: U.F.R. de Math\'{e}matiques, Universit\'{e} Paris 7,
B\^{a}timent Sophie Germain, Case 7012, 75205 Paris Cedex 13, France;
oger@math.univ-paris-diderot.fr.\bigskip 

\noindent \textbf{Abstract.} Generalizing previous results, we give an
algebraic characterization of elementary equivalence for
polycyclic-by-finite groups. We use this characterization to investigate the
relations between their elementary equivalence and the elementary
equivalence of the factors in their decompositions in direct products of
indecomposable groups. In particular we prove that the elementary
equivalence $G\equiv H$ of two such groups $G,H$ is equivalent to each of
the following properties: 1)$\mathrm{~}G\times \cdots \times G\mathrm{\ (}k~%
\mathrm{times~}G\mathrm{)}\equiv H\times \cdots \times H\mathrm{\ (}k~%
\mathrm{times~}H\mathrm{)}$ for an integer $k\geq 1$; 2)$\mathrm{~}A\times
G\equiv B\times H$ for two polycyclic-by-finite groups $A,B$\ such that $%
A\equiv B$. It is not presently known if 1) implies $G\equiv H$ for any
groups $G,H$.\bigskip

\textbf{1. Characterization of elementary equivalence.}\bigskip

In the present paper, we investigate the elementary equivalence between
finitely generated groups and the relations between direct products and
elementary equivalence for groups. First we give the relevant definitions.

We say that two groups $M,N$ are \emph{elementarily equivalent} and we write 
$M\equiv N$ if they satisfy the same first-order sentences in the language
which consists of one binary functional symbol.

We say that a group $M$ is \emph{polycyclic} if there exist some subgroups $%
\left\{ 1\right\} =M_{0}\subset \cdots \subset M_{n}=M$\ with $M_{i-1}$
normal in $M_{i}$ and $M_{i}/M_{i-1}$ cyclic for $1\leq i\leq n$. For any
properties $P,Q$, we say that $M$ is $P$\emph{-by-}$Q$ if there exists a
normal subgroup $N$ which satisfies $P$ and such that $M/N$ satisfies $Q$.

For each group $M$, we denote by $Z(M)$ the center of $M$. For each $%
A\subset M$, we denote by\ $\left\langle A\right\rangle $\ the subgroup of $%
M $ generated by $A$. We write

\noindent $M^{\prime }=\left\langle \left\{ [x,y]\mid x,y\in M\right\}
\right\rangle $\ and, for any $h,k\in 
%TCIMACRO{\U{2115} }%
%BeginExpansion
\mathbb{N}
%EndExpansion
^{\ast }$,

\noindent $M^{\prime }(h)=\{[x_{1},y_{1}]\cdots \lbrack x_{h},y_{h}]\mid
x_{1},y_{1},\ldots ,x_{h},y_{h}\in M\}$,

\noindent $\times ^{k}M=M\times \cdots \times M$ ($k$ times $M)$,

\noindent $M^{k}=\left\langle \left\{ x^{k}\mid x\in M\right\} \right\rangle 
$ and

\noindent $M^{k}(h)=\{x_{1}^{k}\cdots x_{h}^{k}\mid x_{1},\ldots ,x_{h}\in
M\}$.

Each of the properties $M^{\prime }=M^{\prime }(h)$ and $M^{k}=M^{k}(h)$ can
be expressed by a first-order sentence. As polycyclic-by-finite groups are
noetherian, it follows from [16, Corollary 2.6.2] that, for each
polycyclic-by-finite group $M$ and each $k\in 
%TCIMACRO{\U{2115} }%
%BeginExpansion
\mathbb{N}
%EndExpansion
^{\ast }$, there exists $h\in 
%TCIMACRO{\U{2115} }%
%BeginExpansion
\mathbb{N}
%EndExpansion
^{\ast }$\ such that $M^{\prime }=M^{\prime }(h)$ and $M^{k}=M^{k}(h)$.
Actually, the result concerning the property $M^{\prime }=M^{\prime }(h)$
was first proved in [15].

For each group $M$ and any subsets $A,B$, we write $AB=\{xy\mid x\in A$ and $%
y\in B\}$. If $A,B$ are subgroups of $M$ and if $A$ or $B$ is normal, then $%
AB$ is also a subgroup.

For each group $M$ and each subgroup $S$ of $M$ such that $M^{\prime
}\subset S$, we consider the \emph{isolator} $I_{M}(S)=\cup _{k\in 
%TCIMACRO{\U{2115} }%
%BeginExpansion
\mathbb{N}
%EndExpansion
^{\ast }}\{x\in M\mid x^{k}\in S\}$, which is a subgroup of $M$. We write $%
\Gamma (M)=I_{M}(Z(M)M^{\prime })$ and $\Delta (M)=I_{M}(M^{\prime })$. We
also have $\Gamma (M)=I_{M}(Z(M)\Delta (M))$. If $M$ is abelian, then $%
\Delta (M)$ is the \emph{torsion subgroup} $\tau (M)$.

For any $h,k\in 
%TCIMACRO{\U{2115} }%
%BeginExpansion
\mathbb{N}
%EndExpansion
^{\ast }$, there exists a set of first-order sentences which, in each group $%
M$ with $M^{\prime }=M^{\prime }(h)$, expresses that $\Gamma (M)^{k}\subset
Z(M)M^{\prime }$ (resp. $\Delta (M)^{k}\subset M^{\prime }$).\bigskip

During the 2000s, there was a lot of progress in the study of elementary
equivalence between finitely generated groups, with the proof of Tarski's
conjecture stating that all free groups with at least two generators are
elementarily equivalent (see [5] and [17]).

Anyway, examples of elementarily equivalent nonisomorphic finitely generated
groups also exist for groups which are much nearer to the abelian class. One
of them was given as early as 1971 by B.I. Zil'ber in [18].

It quickly appeared that the class of polycyclic-by-finite groups would be
an appropriate setting for an algebraic characterization of elementary
equivalence. Actually, two elementarily equivalent such groups necessarily
have the same finite images (see [9, Remark, p. 475]) and, by [4], any class
of such groups which have the same finite images is a finite union of
isomorphism classes.

Between 1981 and 1991, algebraic characterizations of elementary equivalence
were given for particular classes of polycyclic-by-finite groups (see [9],
[11] and [12]). During the 2000s, some new progress was made with the
introduction of quasi-finitely axiomatizable finitely generated groups by A.
Nies (see [6, Introduction]). This line of investigation, and in particular
the results of [6], made it possible to obtain the characterization of
elementary equivalence for polycyclic-by-finite groups which is given
below.\bigskip

By Feferman-Vaught's theorem (see [3]), we have $G_{1}\times G_{2}\equiv
H_{1}\times H_{2}$ for any groups $G_{1},G_{2},H_{1},H_{2}$ such that $%
G_{1}\equiv H_{1}$ and $G_{2}\equiv H_{2}$. In particular, we have $\times
^{k}G\equiv \times ^{k}H$ for each $k\in 
%TCIMACRO{\U{2115} }%
%BeginExpansion
\mathbb{N}
%EndExpansion
$ and any groups $G,H$ such that $G\equiv H$.

In [7, p. 9], L. Manevitz proposes the following conjecture:\bigskip

\noindent \textbf{Conjecture.} For any groups $G,H$ and each integer $k\geq
2 $, if $\times ^{k}G\equiv \times ^{k}H$, then $G\equiv H$.\bigskip

L. Manevitz mentions that the conjecture becomes false if we consider sets
with one unary function instead of groups. He obtains a counterexample by
taking for $G$ the set $%
%TCIMACRO{\U{2115} }%
%BeginExpansion
\mathbb{N}
%EndExpansion
$\ equipped with the successor function and for $H$\ the disjoint union of
two copies of $G$. We have $\times ^{k}G\cong \times ^{k}H$ for each integer 
$k\geq 2$, but $G$ and $H$ are not elementarily equivalent.

It seems likely that the conjecture is also false for groups in general.
Anyway, it appears that no counterexample has been given until present.

For some classes of groups, we prove that the conjecture is true by showing
that the groups are characterized up to elementary equivalence by invariants
with good behaviour relative to products. We show in this way that the
conjecture is true for abelian groups, using the invariants introduced by
Szmielew or Eklof-Fisher (see [2]).

Similarly, by [11, Cor., p. 1042], two finitely generated abelian-by-finite
groups $G,H$ are elementarily equivalent if and only if they have the same
finite images. For polycyclic-by-finite groups, and in particular for
finitely generated abelian-by-finite groups, the last property is true if
and only if $G/G^{n}\cong H/H^{n}$ for each $n\in 
%TCIMACRO{\U{2115} }%
%BeginExpansion
\mathbb{N}
%EndExpansion
^{\ast }$. The conjecture is true for finitely generated abelian-by-finite
groups because, for any such groups $G,H$ and any $k,n\in 
%TCIMACRO{\U{2115} }%
%BeginExpansion
\mathbb{N}
%EndExpansion
^{\ast }$, the finite groups $G/G^{n}$ and $H/H^{n}$ are isomorphic if and
only if $\times ^{k}(G/G^{n})\cong (\times ^{k}G)/(\times ^{k}G)^{n}$ and $%
\times ^{k}(H/H^{n})\cong (\times ^{k}H)/(\times ^{k}H)^{n}$ are isomorphic.

Also, according to [12, Th., p. 173], which answered a conjecture stated in
[8], two finitely generated finite-by-nilpotent groups $G,H$ are
elementarily equivalent if and only if $G\times 
%TCIMACRO{\U{2124} }%
%BeginExpansion
\mathbb{Z}
%EndExpansion
$ and $H\times 
%TCIMACRO{\U{2124} }%
%BeginExpansion
\mathbb{Z}
%EndExpansion
$ are isomorphic. It follows that Manevitz's conjecture is true for finitely
generated finite-by-nilpotent groups (see [12, Cor. 1, p. 180]).

In the more general case of polycyclic-by-finite groups, neither of the
properties above is a characterisation of elementary equivalence. Actually
(see [12, p. 173]), there exist:

\noindent 1) two finitely generated torsion-free nilpotent groups of class $%
2 $, $G,H$, which have the same finite images and do not satisfy $G\times 
%TCIMACRO{\U{2124} }%
%BeginExpansion
\mathbb{Z}
%EndExpansion
\cong H\times 
%TCIMACRO{\U{2124} }%
%BeginExpansion
\mathbb{Z}
%EndExpansion
$ (for such groups, $G\times 
%TCIMACRO{\U{2124} }%
%BeginExpansion
\mathbb{Z}
%EndExpansion
\cong H\times 
%TCIMACRO{\U{2124} }%
%BeginExpansion
\mathbb{Z}
%EndExpansion
$ implies $G\cong H$);

\noindent 2) two polycyclic abelian-by-finite groups $G,H$ which have the
same finite images and do not satisfy $G\times 
%TCIMACRO{\U{2124} }%
%BeginExpansion
\mathbb{Z}
%EndExpansion
\cong H\times 
%TCIMACRO{\U{2124} }%
%BeginExpansion
\mathbb{Z}
%EndExpansion
$.

\noindent It follows from [11, Cor., p. 1042] that $G$ and $H$\ are
elementarily equivalent in the second case and from [12, Th., p. 173] that
they are not elementarily equivalent in the first case.

However, Theorem 1.1 below gives an algebraic characterization of elementary
equivalence for polycyclic-by-finite groups. In the second section of the
present paper, we use it in order to prove that the conjecture is true for
such groups.\bigskip

\noindent \textbf{Theorem 1.1.} Two polycyclic-by-finite groups $G,H$ are
elementarily equivalent if and only if there exist $n\in 
%TCIMACRO{\U{2115} }%
%BeginExpansion
\mathbb{N}
%EndExpansion
^{\ast }$ such that $\left| G/G^{n}\right| =\left| H/H^{n}\right| $ and, for
each $k\in 
%TCIMACRO{\U{2115} }%
%BeginExpansion
\mathbb{N}
%EndExpansion
^{\ast }$, an injective homomorphism $f_{k}:G\rightarrow H$ such that $%
H=f_{k}(G)Z(H^{n})$ and $\left| H:f_{k}(G)\right| $\ is prime to $k$.\bigskip

The following Proposition implies that the case $n=1$ of the
characterization above is equivalent to the stronger property $G\times 
%TCIMACRO{\U{2124} }%
%BeginExpansion
\mathbb{Z}
%EndExpansion
\cong H\times 
%TCIMACRO{\U{2124} }%
%BeginExpansion
\mathbb{Z}
%EndExpansion
$. We shall give the proof of the Proposition first, since it is much
shorter than the proof of the Theorem.\bigskip

\noindent \textbf{Proposition 1.2.} For any groups $G,H$ with $G/G^{\prime }$%
\ and $H/H^{\prime }$\ finitely generated, the following properties are
equivalent:

\noindent 1) $G\times 
%TCIMACRO{\U{2124} }%
%BeginExpansion
\mathbb{Z}
%EndExpansion
\cong H\times 
%TCIMACRO{\U{2124} }%
%BeginExpansion
\mathbb{Z}
%EndExpansion
$;

\noindent 2) for each $k\in 
%TCIMACRO{\U{2115} }%
%BeginExpansion
\mathbb{N}
%EndExpansion
^{\ast }$, there exists an injective homomorphism $f_{k}:G\rightarrow H$
with $H=f_{k}(G)Z(H)$ and $\left| H/f_{k}(G)\right| $ prime to $k$;

\noindent 3) the same property is true for $k=\left| \Gamma (H)/(Z(H)\Delta
(H))\right| .\left| \Delta (H)/H^{\prime }\right| $.\bigskip

\noindent \textbf{Proof of Proposition 1.2.} First we show that 1) implies
2). We consider an isomorphism $f:M\rightarrow N$ with $M=G\times \langle
u\rangle $, $N=H\times \langle w\rangle $\ and $\left\langle u\right\rangle $%
, $\left\langle w\right\rangle $\ infinite. The restriction of $f$ to $G\cap
f^{-1}(H)$\ is an isomorphism from $G\cap f^{-1}(H)$\ to $f(G)\cap H$.

Suppose $f(G)\neq H$. Then $G/(G\cap f^{-1}(H))\cong \left\langle
G,f^{-1}(H)\right\rangle /f^{-1}(H)\subset M/f^{-1}(H)\cong N/H$\ is
infinite cyclic. The same property is true for $H/(f(G)\cap H)$. We consider 
$v\in G$\ such that $G=\left\langle v,G\cap f^{-1}(H)\right\rangle $\ and $%
x\in H$\ such that $H=\left\langle x,f(G)\cap H\right\rangle $.

We write $f(u)=w^{k}x^{l}y$\ and $f(v)=w^{m}x^{n}z$\ with $k,l,m,n\in 
%TCIMACRO{\U{2124} }%
%BeginExpansion
\mathbb{Z}
%EndExpansion
$\ and $y,z\in f(G)\cap H$. We have $u\in Z(M)$, whence $f(u)\in Z(N)$ and $%
x^{l}y\in Z(H)$. The integers $l,n$\ are prime to each other since $f$\
induces an isomorphism from $M/(G\cap f^{-1}(H))$\ to $N/(f(G)\cap H)$.

For each $r\in 
%TCIMACRO{\U{2124} }%
%BeginExpansion
\mathbb{Z}
%EndExpansion
$, we consider the homomorphism $f_{r}:G\rightarrow H$ defined by $%
f_{r}(v)=(x^{l}y)^{r}x^{n}z$\ and $f_{r}(t)=f(t)$\ for $t\in G\cap f^{-1}(H)$%
. We have $f_{r}(t)\in f(t)Z(H)$ for each $t\in G$. It follows $%
N=f(G)Z(N)=f_{r}(G)Z(N)$\ and $H=f_{r}(G)Z(H)$. We also have $f_{r}(v)\in
x^{lr+n}(f(G)\cap H)$ since $f(G)\cap H$ is normal in $H$. Consequently, $%
f_{r}$ is injective.

For each $s\in 
%TCIMACRO{\U{2115} }%
%BeginExpansion
\mathbb{N}
%EndExpansion
^{\ast }$, we can choose $r$\ in such a way that $lr+n$\ is prime to $s$,
since $l$ and $n$\ are prime to each other. Then $\left| H/f_{r}(G)\right| $
is prime to $s$.

Now it suffices to prove that 3) implies 1). The property $H=f_{k}(G)Z(H)$
implies $H^{\prime }=f_{k}(G)^{\prime }=f_{k}(G^{\prime })$. Consequently, $%
\left| \Delta (H)/f_{k}(\Delta (G))\right| $ both divides $\left| \Delta
(H)/H^{\prime }\right| $\ since $H^{\prime }\subset f_{k}(\Delta (G))$\ and $%
\left| H/f_{k}(G)\right| $\ since $f_{k}(\Delta (G))=\Delta (H)\cap f_{k}(G)$%
. As $\left| H/f_{k}(G)\right| $ is prime to $\left| \Delta (H)/H^{\prime
}\right| $, it follows $f_{k}(\Delta (G))=\Delta (H)$.

Now the property $G\times 
%TCIMACRO{\U{2124} }%
%BeginExpansion
\mathbb{Z}
%EndExpansion
\cong H\times 
%TCIMACRO{\U{2124} }%
%BeginExpansion
\mathbb{Z}
%EndExpansion
$ follows from [12, Prop. 3] and its proof. Actually, the group $G$\
considered in that proposition is supposed finitely generated
finite-by-nilpotent, but the proof only uses the property ``$G/\Delta (G)$\
finitely generated''.~~$\blacksquare $\bigskip

The definitions and the two lemmas below will be used in the proof of
Theorem 1.1.

For each group $M$, we write $E(M)=\cap _{k\in 
%TCIMACRO{\U{2115} }%
%BeginExpansion
\mathbb{N}
%EndExpansion
^{\ast }}M^{k}$. If $M$ is abelian and if $\tau (M)$ is finite, then $E(M)$
is the additive structure of a $%
%TCIMACRO{\U{211a} }%
%BeginExpansion
\mathbb{Q}
%EndExpansion
$-vector space.\bigskip

\noindent \textbf{Definitions.} Let $M$ be a group, let $K$ be a group of
automorphisms of $M$ and let $A,B$\ be subgroups of $M$\ stabilized by $K$
such that $A\cap B=\left[ A,B\right] =\left\{ 1\right\} $. We say that $B$
is a $K$\emph{-complement} of $A$\ in $M$ if $AB=M$, and a $K$\emph{%
-quasi-complement} of $A$\ in $M$ if $\left| M:AB\right| $ is finite.\bigskip

\noindent \textbf{Lemma 1.3.} Let $M$ be a finitely generated abelian group
with additive notation, let $N,S$\ be subgroups of $M$ such that $N\cap
S=\left\{ 0\right\} $, and let $K$ be a finite group of automorphisms of $M$%
\ which stabilizes $N$ and $S$. Then there exists a $K$-quasi-complement $T$
of $N$ in $M$\ which contains $S$.\bigskip

\noindent \textbf{Proof of Lemma 1.3.} We consider an integer $r\geq 1$ such
that $r.\tau (M)=\left\{ 0\right\} $ and $r.I_{M}(S)\subset S$. We write $%
M^{\ast }=rM$, $N^{\ast }=M^{\ast }\cap N$\ and $S^{\ast }=M^{\ast }\cap S$.
Then $M^{\ast }$\ and $M^{\ast }/S^{\ast }$\ are torsion-free. We have $%
N^{\ast }\cap S^{\ast }=\left\{ 0\right\} $ and $M^{\ast },N^{\ast },S^{\ast
}$ are stabilized by $K$. As $M/M^{\ast }$\ is finite, it suffices to show
that $N^{\ast }$\ has a $K$-quasi-complement in $M^{\ast }$ which contains $%
S^{\ast }$, or equivalently that $N^{\ast }$\ has a $K$-quasi-complement in $%
M^{\ast }/S^{\ast }$.

Consequently, it suffices to prove the Lemma for $M$ finitely generated
torsion-free abelian and $S=\left\{ 0\right\} $. Then $M$\ can be embedded
in a $%
%TCIMACRO{\U{211a} }%
%BeginExpansion
\mathbb{Q}
%EndExpansion
$-vector space $\overline{M}$\ of finite dimension such that the action of $%
K $ on $\overline{M}$\ is defined. Moreover, $K$\ stabilizes the isolator $%
\overline{N}$ of $N$\ in $\overline{M}$, which is a subspace of $\overline{M}
$. It follows from the proof of Maschke's theorem given in [1, pp. 226-227]
that $\overline{N}$\ has a $K$-complement $T$\ in $\overline{M}$. Then $%
T\cap M$ is a $K$-quasi-complement of $N$ in $M$.~~$\blacksquare $\bigskip

\noindent \textbf{Lemma 1.4.} Let $M$ be a torsion-free abelian group with
additive notation, let $N$ be a divisible subgroup of $M$, let $S$ be a
subgroup of $M$\ such that $N\cap S=\left\{ 0\right\} $, and let $K$ be a
finite group of automorphisms of $M$\ which stabilizes $N$ and $S$. Then
there exists a $K$-complement $T$ of $N$ in $M$\ which contains $S$.\bigskip

\noindent \textbf{Proof of Lemma 1.4.} It suffices to show that $M=N\oplus S$%
\ if $S$ is maximal for the properties $N\cap S=\left\{ 0\right\} $ and $S$\
stabilized by $K$. Then we have $I_{M}(S)=S$ since $N\cap I_{M}(S)=\left\{
0\right\} $ and $I_{M}(S)$\ is stabilized by $K$. We reduce the proof to the
case $S=\left\{ 0\right\} $\ by considering $M/S$ instead of $M$.

Now suppose $S=\left\{ 0\right\} $\ and consider $x\in M-N$. Then $K$
stabilizes $T=\langle \{f(x)\mid f\in K\}\rangle $. In order to obtain a
contradiction, it suffices to prove that there exists a nontrivial subgroup $%
U$ of $T$\ stabilized by $K$\ such that $N\cap U=\left\{ 0\right\} $.\ But,
as $T$ is finitely generated, it follows from Lemma 1.3 that $N\cap T$ has a 
$K$-quasi-complement $U$ in $T$. Moreover, $U$ is nontrivial since $T/(N\cap
T)$ is torsion-free and therefore infinite.~~$\blacksquare $\bigskip

\noindent \textbf{Proof of Theorem 1.1.} First we show that the condition is
necessary. We consider a finite generating tuple $\overline{x}$ of $G$. By
[6, Th. 3.1], there exist an integer $n\geq 1$ and, for each integer $k\geq
1 $, a formula $\varphi _{k}(\overline{u})$ satisfied by $\overline{x}$\ in $%
G$ such that, for each finitely generated group $M$ and each tuple $%
\overline{y} $\ which satisfies $\varphi _{k}$ in $M$, the map $\overline{x}%
\rightarrow \overline{y}$ induces an isomorphism from $G/Z(G^{n})$ to $%
M/Z(M^{n})$ and an injective homomorphism $h:G\rightarrow M$ with $\left|
M:h(G)\right| $\ prime to $k$. As $G$ and $H$ are elementarily equivalent,
we have $\left| G/G^{n}\right| =\left| H/H^{n}\right| $ and there exists a
tuple $\overline{y}$\ which satisfies $\varphi _{k}$ in $H$.

It remains to be proved that the condition is sufficient. First we show some
properties which are true for each injective homomorphism $f:G\rightarrow H$%
\ such that $H=f(G)Z(H^{n})$.

The homomorphism $f$ induces an isomorphism from $G/G^{n}$\ to $H/H^{n}$
since $H=f(G)H^{n}$ and $\left| G/G^{n}\right| =\left| H/H^{n}\right| $.
Consequently, we have $f^{-1}(H^{n})=G^{n}$, $f(G^{n})=H^{n}\cap f(G)$ and $%
H^{n}=(H^{n}\cap f(G))Z(H^{n})=f(G^{n})Z(H^{n})$. It follows $%
(H^{n})^{\prime }=(f(G^{n}))^{\prime }=f((G^{n})^{\prime })$ and $%
f^{-1}(\Delta (H^{n}))=\Delta (G^{n})$.

Now we show that $f^{-1}(Z(H^{n}))=Z(G^{n})$.\ For each $x\in G$\ such that $%
f(x)\in Z(H^{n})$, we have $x\in G^{n}$\ since $f(x)\in H^{n}$, and there
exists no $y\in G^{n}$\ such that $\left[ x,y\right] \neq 1$\ since it would
imply $f(y)\in H^{n}$ and $\left[ f(x),f(y)\right] \neq 1$. Conversely, for
each $x\in Z(G^{n})$, we have $f(x)\in H^{n}$. Let us suppose that there
exists $z\in H^{n}$\ such that $\left[ f(x),z\right] \neq 1$. Then we have $%
z=z_{1}^{n}\cdots z_{r}^{n}$ with $r\in 
%TCIMACRO{\U{2115} }%
%BeginExpansion
\mathbb{N}
%EndExpansion
^{\ast }$ and $z_{1},\ldots ,z_{r}\in H$. For each $i\in \left\{ 1,\ldots
,r\right\} $, we consider $y_{i}\in G$ such that $z_{i}\in f(y_{i})Z(H^{n})$%
. As $Z(H^{n})$ is normal in $H$, we have $z\in f(y_{1})^{n}\cdots
f(y_{r})^{n}Z(H^{n})$, whence $f(\left[ x,y_{1}^{n}\cdots y_{r}^{n}\right] )=%
\left[ f(x),f(y_{1})^{n}\cdots f(y_{r})^{n}\right] =\left[ f(x),z\right]
\neq 1$ and $\left[ x,y_{1}^{n}\cdots y_{r}^{n}\right] \neq 1$, which
contradicts $x\in Z(G^{n})$.

Consequently, we have $f(Z(G^{n}))=Z(H^{n})\cap f(G)$\ and $f$ induces an
isomorphism from $G/Z(G^{n})$ to $H/Z(H^{n})$. It follows $\left|
H:f(G)\right| =\left| Z(H^{n})/f(Z(G^{n}))\right| $.

Now suppose that $\left| H:f(G)\right| $ is prime to $\left| \Delta
(H^{n})/(H^{n})^{\prime }\right| $. Then $f^{-1}(\Delta (H^{n}))=\Delta
(G^{n})$\ implies $f(\Delta (G^{n}))=\Delta (H^{n})$ since $\left| \Delta
(H^{n})/f(\Delta (G^{n}))\right| $\ both divides $\left| \Delta
(H^{n})/(H^{n})^{\prime }\right| $ and $\left| H:f(G)\right| $. Moreover,
the equalities $f(Z(G^{n}))=Z(H^{n})\cap f(G)$ and $f(\Delta (G^{n}))=\Delta
(H^{n})$ imply $f(Z(G^{n})\cap \Delta (G^{n}))=Z(H^{n})\cap \Delta (H^{n})$.

For each $k\in 
%TCIMACRO{\U{2115} }%
%BeginExpansion
\mathbb{N}
%EndExpansion
$, we consider an injective homomorphism $f_{k}:G\rightarrow H$ with $%
H=f_{k}(G)Z(H^{n})$ and $\left| H:f_{k}(G)\right| $ prime to $l!$ where $%
l=\sup (k,\left| \Delta (H^{n})/(H^{n})^{\prime }\right| )$. Each $f_{k}$
satisfies $f_{k}^{-1}(H^{n})=G^{n}$, $f_{k}^{-1}(Z(H^{n}))=Z(G^{n})$, $%
H^{n}=f_{k}(G^{n})Z(H^{n})$, $f_{k}(\Delta (G^{n}))=\Delta (H^{n})$, $%
f_{k}(Z(G^{n})\cap \Delta (G^{n}))=Z(H^{n})\cap \Delta (H^{n})$ and $\left|
H:f_{k}(G)\right| =\left| Z(H^{n})/f_{k}(Z(G^{n}))\right| $.

We consider an $\omega _{1}$-incomplete ultrafilter\ $\mathcal{U}$\ over $%
%TCIMACRO{\U{2115} }%
%BeginExpansion
\mathbb{N}
%EndExpansion
$ and the injective homomorphism $f:G^{\mathcal{U}}\rightarrow H^{\mathcal{U}%
}$\ which admits $(f_{k})_{k\in 
%TCIMACRO{\U{2115} }%
%BeginExpansion
\mathbb{N}
%EndExpansion
}$ as a representative.

We have $(G^{\mathcal{U}})^{n}=(G^{n})^{\mathcal{U}}$ and $(H^{\mathcal{U}%
})^{n}=(H^{n})^{\mathcal{U}}$ since there exists an integer $r\geq 1$ such
that $G^{n}=G^{n}(r)$\ and $H^{n}=H^{n}(r)$. Consequently, $f$\ induces an
injective homomorphism from $(G^{\mathcal{U}})^{n}$ to $(H^{\mathcal{U}%
})^{n} $ and an isomorphism from $G^{\mathcal{U}}/(G^{\mathcal{U}})^{n}$ to $%
H^{\mathcal{U}}/(H^{\mathcal{U}})^{n}$.

We also have $((G^{\mathcal{U}})^{n})^{\prime }=((G^{n})^{\mathcal{U}%
})^{\prime }=((G^{n})^{\prime })^{\mathcal{U}}$ and $((H^{\mathcal{U}%
})^{n})^{\prime }=((H^{n})^{\mathcal{U}})^{\prime }=((H^{n})^{\prime })^{%
\mathcal{U}}$ since there exists an integer $s\geq 1$ such that $%
(G^{n})^{\prime }=(G^{n})^{\prime }(s)$ and $(H^{n})^{\prime
}=(H^{n})^{\prime }(s)$. It follows $\Delta ((G^{\mathcal{U}})^{n})=\Delta
(G^{n})^{\mathcal{U}}$ and $\Delta ((H^{\mathcal{U}})^{n})=\Delta (H^{n})^{%
\mathcal{U}}$ since there exists an integer $t\geq 1$ such that $\Delta
(G^{n})^{t}\subset (G^{n})^{\prime }$ and $\Delta (H^{n})^{t}\subset
(H^{n})^{\prime }$. Now the equalities $f_{k}(\Delta (G^{n}))=\Delta (H^{n})$
for $k\in 
%TCIMACRO{\U{2115} }%
%BeginExpansion
\mathbb{N}
%EndExpansion
$\ imply $f(\Delta ((G^{\mathcal{U}})^{n}))=\Delta ((H^{\mathcal{U}})^{n})$.

We have $Z((G^{\mathcal{U}})^{n})=Z((G^{n})^{\mathcal{U}})=Z(G^{n})^{%
\mathcal{U}}$ and $Z((H^{\mathcal{U}})^{n})=Z(H^{n})^{\mathcal{U}}$. The
equalities $f_{k}(Z(G^{n})\cap \Delta (G^{n}))=Z(H^{n})\cap \Delta (H^{n})$
imply $f(Z((G^{\mathcal{U}})^{n})\cap \Delta ((G^{\mathcal{U}})^{n}))=Z((H^{%
\mathcal{U}})^{n})\cap \Delta ((H^{\mathcal{U}})^{n})$.

The action of $G/G^{n}$ on $Z(G^{n})$ is defined since $G^{n}$ acts
trivially on $Z(G^{n})$, and the action of $G^{\mathcal{U}}/(G^{\mathcal{U}%
})^{n}$ on $Z((G^{\mathcal{U}})^{n})$ is defined since $(G^{\mathcal{U}%
})^{n} $ acts trivially on $Z((G^{\mathcal{U}})^{n})$.\ The same properties
are true for $H$.

According to Lemma 1.3, there exist a $G/G^{n}$-quasi-complement $C$ of $%
Z(G^{n})\cap \Delta (G^{n})$\ in $Z(G^{n})$ and an $H/H^{n}$%
-quasi-complement $D$ of $Z(H^{n})\cap \Delta (H^{n})$\ in $Z(H^{n})$.

Then $C^{\mathcal{U}}$ is a $G^{\mathcal{U}}/(G^{\mathcal{U}})^{n}$%
-quasi-complement of $Z((G^{\mathcal{U}})^{n})\cap \Delta ((G^{\mathcal{U}%
})^{n})$ in $Z((G^{\mathcal{U}})^{n})$, and $E(C^{\mathcal{U}})$ is a $G^{%
\mathcal{U}}/(G^{\mathcal{U}})^{n}$-complement of $E(Z((G^{\mathcal{U}%
})^{n})\cap \Delta ((G^{\mathcal{U}})^{n}))$ in $E(Z((G^{\mathcal{U}})^{n}))$
since $E(Z((G^{\mathcal{U}})^{n}))$ is a $%
%TCIMACRO{\U{211a} }%
%BeginExpansion
\mathbb{Q}
%EndExpansion
$-vector space. In particular, $E(C^{\mathcal{U}})$ is a normal subgroup of $%
G^{\mathcal{U}}$. The same properties are true for $H$ and $D$.

Now we prove that there exists a subgroup $T$ of $G^{\mathcal{U}}$\ such
that $G\subset T$, $T\cap E(C^{\mathcal{U}})=\left\{ 1\right\} $ and $G^{%
\mathcal{U}}=T.E(C^{\mathcal{U}})$. As $E(C^{\mathcal{U}})\cap \Delta ((G^{%
\mathcal{U}})^{n})=\left\{ 1\right\} $, it suffices to show that there
exists a subgroup $S$ of\ $G^{\mathcal{U}}/\Delta ((G^{\mathcal{U}})^{n})$\
such that $G/\Delta (G)\subset S$, $S\cap E(C^{\mathcal{U}})=\left\{
1\right\} $\ and $S.E(C^{\mathcal{U}})=G^{\mathcal{U}}/\Delta ((G^{\mathcal{U%
}})^{n})$.

The action of $G^{\mathcal{U}}/(G^{\mathcal{U}})^{n}$ on $(G^{\mathcal{U}%
})^{n}/\Delta ((G^{\mathcal{U}})^{n})$\ is defined since $(G^{\mathcal{U}%
})^{n}/\Delta ((G^{\mathcal{U}})^{n})$\ is an abelian normal subgroup of $G^{%
\mathcal{U}}/\Delta ((G^{\mathcal{U}})^{n})$.\ By Lemma 1.4, $E(C^{\mathcal{U%
}})$ has a $G^{\mathcal{U}}/(G^{\mathcal{U}})^{n}$-complement $R$ in $(G^{%
\mathcal{U}})^{n}/\Delta ((G^{\mathcal{U}})^{n})$\ which contains $%
G^{n}/\Delta (G^{n})$.

Let us consider $x_{1},\ldots ,x_{m}\in G$ such that $G$\ is the disjoint
union of $x_{1}G^{n},\ldots ,x_{m}G^{n}$. Then $S=x_{1}R\cup \ldots \cup
x_{m}R$ is a subgroup of $G^{\mathcal{U}}/\Delta ((G^{\mathcal{U}})^{n})$\
which contains $G/\Delta (G^{n})$. We have $S\cap E(C^{\mathcal{U}})=\left\{
1\right\} $\ and $S.E(C^{\mathcal{U}})=G^{\mathcal{U}}/\Delta ((G^{\mathcal{U%
}})^{n})$.

For each $k\in 
%TCIMACRO{\U{2115} }%
%BeginExpansion
\mathbb{N}
%EndExpansion
$, $f_{k}$ induces an isomorphism from $Z(G^{n})/Z(G^{n})^{k!}$ to $%
Z(H^{n})/Z(H^{n})^{k!}$\ since $\left| Z(H^{n})/f_{k}(Z(G^{n}))\right|
=\left| H:f_{k}(G)\right| $ is prime to $k!$. Consequently, $f$ induces an
isomorphism from $Z((G^{\mathcal{U}})^{n})/E(Z((G^{\mathcal{U}})^{n}))$ to

\noindent $Z((H^{\mathcal{U}})^{n})/E(Z((H^{\mathcal{U}})^{n}))$, and
therefore induces an isomorphism from

\noindent $G^{\mathcal{U}}/E(Z((G^{\mathcal{U}})^{n}))$ to $H^{\mathcal{U}%
}/E(Z((H^{\mathcal{U}})^{n}))$ since it induces an isomorphism from $G^{%
\mathcal{U}}/Z((G^{\mathcal{U}})^{n})$ to $H^{\mathcal{U}}/Z((H^{\mathcal{U}%
})^{n})$.

It follows that $f$ induces an isomorphism from $G^{\mathcal{U}}/E(C^{%
\mathcal{U}})$ to $H^{\mathcal{U}}/E(D^{\mathcal{U}})$\ since\ $f(E(Z((G^{%
\mathcal{U}})^{n})\cap \Delta ((G^{\mathcal{U}})^{n})))=E(Z((H^{\mathcal{U}%
})^{n})\cap \Delta ((H^{\mathcal{U}})^{n}))$. In particular, each element of 
$H^{\mathcal{U}}$ can be written in a unique way as a product of an element
of $f(T)$ and an element of $E(D^{\mathcal{U}})$.

Now denote by $u_{1},\ldots ,u_{m}$\ the elements of $G^{\mathcal{U}}/(G^{%
\mathcal{U}})^{n}$ and write $v_{i}=f(u_{i})$\ for $1\leq i\leq m$. The
restriction of $f$ to $T$ is completed into an isomorphism from $G^{\mathcal{%
U}}$ to $H^{\mathcal{U}}$ by any isomorphism $g:E(C^{\mathcal{U}%
})\rightarrow E(D^{\mathcal{U}})$ which satisfies $g(z^{u_{i}})=g(z)^{v_{i}}$
for $z\in E(C^{\mathcal{U}})$\ and $1\leq i\leq m$. It remains to be proved
that $g$ exists.

The automorphisms $\theta _{u_{i}}:x\rightarrow x^{u_{i}}$ of $%
Z(G^{n})/(Z(G^{n})\cap \Delta (G^{n}))$ and the automorphisms $\theta
_{v_{i}}:x\rightarrow x^{v_{i}}$ of $Z(H^{n})/(Z(H^{n})\cap \Delta (H^{n}))$
are defined since $G^{n}$ acts trivially on $Z(G^{n})$ and $H^{n}$ acts
trivially on $Z(H^{n})$.

Each $f_{k}$ induces an injective homomorphism from

\noindent $A=(Z(G^{n})/(Z(G^{n})\cap \Delta (G^{n})),\theta _{u_{1}},\ldots
,\theta _{u_{m}})$ to

\noindent $B=(Z(H^{n})/(Z(H^{n})\cap \Delta (H^{n})),\theta _{v_{1}},\ldots
,\theta _{v_{m}})$

\noindent and therefore induces an isomorphism from $A/A^{k!}$ to $B/B^{k!}$
since

\noindent $\left| Z(H^{n})/f_{k}(Z(G^{n}))\right| =\left| H:f_{k}(G)\right| $
is prime to $k!$.

By [13, Cor. 1.2], $A$ and $B$ are elementarily equivalent. Moreover, the
existence of an isomorphism from $A^{\mathcal{U}}$ to $B^{\mathcal{U}}$
follows from the proofs of [11] and [13], and from the existence of
isomorphisms $\widehat{%
%TCIMACRO{\U{2124} }%
%BeginExpansion
\mathbb{Z}
%EndExpansion
}^{\mathcal{U}}\rightarrow 
%TCIMACRO{\U{2124} }%
%BeginExpansion
\mathbb{Z}
%EndExpansion
^{\mathcal{U}}\rightarrow \widehat{%
%TCIMACRO{\U{2124} }%
%BeginExpansion
\mathbb{Z}
%EndExpansion
}\times E(%
%TCIMACRO{\U{2124} }%
%BeginExpansion
\mathbb{Z}
%EndExpansion
^{\mathcal{U}})$ which fix the elements of $%
%TCIMACRO{\U{2124} }%
%BeginExpansion
\mathbb{Z}
%EndExpansion
$, where $\widehat{%
%TCIMACRO{\U{2124} }%
%BeginExpansion
\mathbb{Z}
%EndExpansion
}$ is the profinite completion of $%
%TCIMACRO{\U{2124} }%
%BeginExpansion
\mathbb{Z}
%EndExpansion
$.

Any isomorphism\ from

\noindent $A^{\mathcal{U}}=(Z((G^{\mathcal{U}})^{n})/(Z((G^{\mathcal{U}%
})^{n})\cap \Delta ((G^{\mathcal{U}})^{n})),\theta _{u_{1}},\ldots ,\theta
_{u_{m}})$ to

\noindent $B^{\mathcal{U}}=(Z((H^{\mathcal{U}})^{n})/(Z((H^{\mathcal{U}%
})^{n})\cap \Delta ((H^{\mathcal{U}})^{n})),\theta _{v_{1}},\ldots ,\theta
_{v_{m}})$

\noindent induces an isomorphism from $E(A^{\mathcal{U}})$ to $E(B^{\mathcal{%
U}})$, and therefore induces an isomorphism $g:E(C^{\mathcal{U}})\rightarrow
E(D^{\mathcal{U}})$ which satisfies $g(z^{u_{i}})=g(z)^{v_{i}}$ for $z\in
E(C^{\mathcal{U}})$\ and $1\leq i\leq m$.~~$\blacksquare $\bigskip

\textbf{2. Direct products and elementary equivalence.}\bigskip

For each group $G$ and any subgroups $A,B$, we write $G=A\times B$\ if $G=AB$%
\ and $A\cap B=[A,B]=\left\{ 1\right\} $.\bigskip

\noindent \textbf{Theorem 2.1.} Consider a polycyclic-by-finite group $G$,
an integer $k\geq 1$ such that $\Gamma (G)^{k}\subset Z(G)\Delta (G)$ and $%
\Delta (G)^{k}\subset G^{\prime }$, an integer $n\geq 1$ and some finite
sequences $\overline{x}_{1},\ldots ,\overline{x}_{n}\subset G$ such that $%
G=\left\langle \overline{x}_{1}\right\rangle \times \cdots \times
\left\langle \overline{x}_{n}\right\rangle $. Write $G_{i}=\left\langle 
\overline{x}_{i}\right\rangle $ for $1\leq i\leq n$.\ Then there exist an
integer $r\geq 1$ and, for each integer $m\geq 1$, a formula $\varphi _{m}(%
\overline{u}_{1},\ldots ,\overline{u}_{n})$ satisfied by $(\overline{x}%
_{1},\ldots ,\overline{x}_{n})$\ in $G$ such that, for each finitely
generated group $H$ with $\Gamma (H)^{k}\subset Z(H)\Delta (H)$ and $\Delta
(H)^{k}\subset H^{\prime }$, and for each $(\overline{y}_{1},\ldots ,%
\overline{y}_{n})$\ which satisfies $\varphi _{m}$\ in $H$, there exist some
subgroups $H_{i}\subset C_{H}(\overline{y}_{1},\ldots ,\overline{y}_{i-1},%
\overline{y}_{i+1},\ldots ,\overline{y}_{n})$ with $\left|
H_{i}/H_{i}^{r}\right| =\left| G_{i}/G_{i}^{r}\right| $\ such that $%
H=H_{1}\times \cdots \times H_{n}$ and such that the maps $\overline{x}%
_{i}\rightarrow \mathrm{pr}_{H_{i}}(\overline{y}_{i})$ induce injective
homomorphisms $f_{i}:G_{i}\rightarrow H_{i}$ with $%
H_{i}=f_{i}(G_{i})Z(H_{i}^{r})$ and $\left| H_{i}:f_{i}(G_{i})\right| $\
prime to $m$.\bigskip

\noindent \textbf{Proof.} By [6, Th. 3.1], there exist an integer $r\geq 1$
and, for each integer $m\geq 1$, a formula $\psi _{m}(\overline{u}%
_{1},\ldots ,\overline{u}_{n})$ satisfied by $(\overline{x}_{1},\ldots ,%
\overline{x}_{n})$ in $G$ such that, for each finitely generated group $H$
and each $(\overline{u}_{1},\ldots ,\overline{u}_{n})$\ which satisfies $%
\psi _{m}$\ in $H$, the map $(\overline{x}_{1},\ldots ,\overline{x}%
_{n})\rightarrow (\overline{u}_{1},\ldots ,\overline{u}_{n})$ induces an
injective homomorphism $f:G\rightarrow H$ with $H=f(G)Z(H^{r})$ and $\left|
H:f(G)\right| $\ prime to $m$.

It suffices to show that $\varphi _{m}$\ exists for $m$ divisible by $r$. We
fix $m$ for the remainder of the proof.

For $1\leq i\leq n$, we write $\overline{x}_{i}^{\ast }=(\overline{x}%
_{1},\ldots ,\overline{x}_{i-1},\overline{x}_{i+1},\ldots ,\overline{x}_{n})$
and $\overline{u}_{i}^{\ast }=(\overline{u}_{1},\ldots ,\overline{u}_{i-1},%
\overline{u}_{i+1},\ldots ,\overline{u}_{n})$. We have

\noindent $G_{i}\subset C_{G}(\overline{x}_{i}^{\ast })=Z(G_{1})\times
\cdots \times Z(G_{i-1})\times G_{i}\times Z(G_{i+1})\times \cdots \times
Z(G_{n})=G_{i}Z(G)$.

\noindent It follows that $G=C_{G}(\overline{x}_{1}^{\ast })\cdots C_{G}(%
\overline{x}_{n}^{\ast })$.

We consider an integer $t\geq 1$ such that $G^{\prime }=G^{\prime }(t)$.\ We
denote by $\alpha (\overline{u}_{1},\ldots ,\overline{u}_{n})$ a formula
which says that $H^{\prime }=H^{\prime }(t)$, $[C_{H}(\overline{u}_{i}^{\ast
}),C_{H}(\overline{u}_{j}^{\ast })]=\left\{ 1\right\} $ for $i\neq j$ and $%
H=C_{H}(\overline{u}_{1}^{\ast })\cdots C_{H}(\overline{u}_{n}^{\ast })$. It
follows that $Z(C_{H}(\overline{u}_{i}^{\ast }))=Z(H)$ for $1\leq i\leq n$
and $H/Z(H)=C_{H}(\overline{u}_{1}^{\ast })/Z(H)\times \cdots \times C_{H}(%
\overline{u}_{n}^{\ast })/Z(H)$. The formula $\alpha $\ is satisfied by $(%
\overline{x}_{1},\ldots ,\overline{x}_{n})$ in $G$.

From now on, we only consider groups $H$ such that $H^{\prime }=H^{\prime
}(t)$, $\Gamma (H)^{k}\subset Z(H)\Delta (H)$ and $\Delta (H)^{k}\subset
H^{\prime }$, and sequences $(\overline{u}_{1},\ldots ,\overline{u}_{n})$
which satisfy $\alpha \wedge \psi _{m}$\ in $H$. For each $i\in \left\{
1,\ldots ,n\right\} $, we consider a sequence of terms $\overline{\rho }_{i}(%
\overline{u}_{i})$ such that each element of $\Delta (G_{i})$\ can be
written in a unique way in the form $xy$ with $x\in \overline{\rho }_{i}(%
\overline{x}_{i})$ and $y\in G_{i}^{\prime }$.

As $\alpha (\overline{u}_{1},\ldots ,\overline{u}_{n})$ implies $C_{H}(%
\overline{u}_{i}^{\ast })^{\prime }=C_{H}(\overline{u}_{i}^{\ast })^{\prime
}(t)$ for $1\leq i\leq n$, there exists a formula $\beta (\overline{u}%
_{1},\ldots ,\overline{u}_{n})$ satisfied by $(\overline{x}_{1},\ldots ,%
\overline{x}_{n})$ in $G$ which expresses that:

\noindent 1) $H^{\prime }=C_{H}(\overline{u}_{1}^{\ast })^{\prime }\times
\cdots \times C_{H}(\overline{u}_{n}^{\ast })^{\prime }$ (it suffices to say
that

\noindent $\left\{ 1\right\} =C_{H}(\overline{u}_{i}^{\ast })^{\prime }\cap
(C_{H}(\overline{u}_{1}^{\ast })^{\prime }\cdots C_{H}(\overline{u}%
_{i-1}^{\ast })^{\prime }C_{H}(\overline{u}_{i+1}^{\ast })^{\prime }\cdots
C_{H}(\overline{u}_{n}^{\ast })^{\prime })$ for $1\leq i\leq n$);

\noindent 2) $vw\in \overline{\rho }_{i}(\overline{u}_{i})C_{H}(\overline{u}%
_{i}^{\ast })^{\prime }$ for $1\leq i\leq n$ and $v,w\in \overline{\rho }%
_{i}(\overline{u}_{i})$;

\noindent 3) each element of $\Delta (H)$\ can be written in a unique way as 
$v_{1}\cdots v_{n}v$ with $v\in H^{\prime }$ and $v_{i}\in \overline{\rho }%
_{i}(\overline{u}_{i})$ for $1\leq i\leq n$ (here we use $\Delta
(H)^{k}\subset H^{\prime }$).

\noindent It follows from 1), 2), 3) that $\Delta (H)$ is the direct product
of the subgroups $\overline{\rho }_{i}(\overline{u}_{i})C_{H}(\overline{u}%
_{i}^{\ast })^{\prime }$ for $1\leq i\leq n$.

From now on, we only consider sequences $(\overline{u}_{1},\ldots ,\overline{%
u}_{n})$ which satisfy $\alpha \wedge \beta \wedge \psi _{m}$\ in $H$.

For each $i\in \left\{ 1,\ldots ,n\right\} $, we consider some terms $\sigma
_{i,1}(\overline{u}_{i}),\ldots ,\sigma _{i,r(i)}(\overline{u}_{i})$ such
that each element of $\Gamma (G_{i})$ can be written in a unique way as

\noindent $\sigma _{i,1}(\overline{x}_{i})^{a_{1}}\cdots \sigma _{i,r(i)}(%
\overline{x}_{i})^{a_{r(i)}}y$ with $a_{1},\ldots ,a_{r(i)}\in 
%TCIMACRO{\U{2124} }%
%BeginExpansion
\mathbb{Z}
%EndExpansion
$ and $y\in \Delta (G_{i})$.

There exists a formula $\gamma (\overline{u}_{1},\ldots ,\overline{u}_{n})$
satisfied by $(\overline{x}_{1},\ldots ,\overline{x}_{n})$ in $G$ which
expresses that each element of $\Gamma (H)$ can be written as $u=\Pi _{1\leq
i\leq n}^{1\leq j\leq r(i)}\sigma _{i,j}(\overline{u}_{i})^{a_{i,j}}v^{km}w$
with $0\leq a_{i,j}\leq km-1$\ completely determined by $u$ for any $i,j$,
and $v\in \Gamma (H)$, $w\in \Delta (H)$.

For each $i\in \left\{ 1,\ldots ,n\right\} $, we also consider some terms $%
\tau _{i,1}(\overline{u}_{i}),\ldots ,\tau _{i,s(i)}(\overline{u}_{i})$ such
that each element of $G_{i}$ can be written in a unique way as

\noindent $\tau _{i,1}(\overline{x}_{i})^{a_{1}}\cdots \tau _{i,s(i)}(%
\overline{x}_{i})^{a_{s(i)}}y$ with $a_{1},\ldots ,a_{s(i)}\in 
%TCIMACRO{\U{2124} }%
%BeginExpansion
\mathbb{Z}
%EndExpansion
$ and $y\in \Gamma (G_{i})$. It follows that each element of $C_{G}(%
\overline{x}_{i}^{\ast })=G_{i}Z(G)$\ can be written as

\noindent $x=\tau _{i,1}(\overline{x}_{i})^{a_{1}}\cdots \tau _{i,s(i)}(%
\overline{x}_{i})^{a_{s(i)}}y^{m}z$ with $0\leq a_{1},\ldots ,a_{s(i)}\leq
m-1$\ completely determined by $x$, and $y\in C_{G}(\overline{x}_{i}^{\ast
}) $, $z\in \Gamma (C_{G}(\overline{x}_{i}^{\ast }))$.

Consequently, there exists a formula $\delta _{i}(\overline{u}_{1},\ldots ,%
\overline{u}_{n})$\ satisfied by $(\overline{x}_{1},\ldots ,\overline{x}%
_{n}) $\ in $G$ which expresses that each element of $C_{H}(\overline{u}%
_{i}^{\ast })$ can be written as $u=\tau _{i,1}(\overline{u}%
_{i})^{a_{1}}\cdots \tau _{i,s(i)}(\overline{u}_{i})^{a_{s(i)}}v^{m}w$ with $%
0\leq a_{1},\ldots ,a_{s(i)}\leq m-1$\ completely determined by $u$, and $%
v\in C_{H}(\overline{u}_{i}^{\ast })$, $w\in \Gamma (C_{H}(\overline{u}%
_{i}^{\ast }))$.

We denote by $\delta $ the conjunction of the formulas $\delta _{i}$. From
now on, we only consider sequences $(\overline{u}_{1},\ldots ,\overline{u}%
_{n})$ which satisfy $\varphi _{m}=\alpha \wedge \beta \wedge \gamma \wedge
\delta \wedge \psi _{m}$\ in $H$.

As the abelian group $\Gamma (H)/\Delta (H)$ is finitely generated and
torsion-free, it is freely generated by the images of a family of elements $%
(v_{i,j})_{1\leq i\leq n}^{1\leq j\leq r(i)}\subset \Gamma (H)$, where each $%
v_{i,j}$ can be written as $\sigma _{i,j}(\overline{u}_{i})v^{km}w$ with $%
v\in \Gamma (H)$ and $w\in \Delta (H)$, which implies $v^{km}w\in
Z(H)^{m}\Delta (H)$ since $\Gamma (H)^{k}\subset Z(H)\Delta (H)$. For any $%
i,j$, we choose $w$ in such a way that $v^{km}w\in Z(H)^{m}$, which implies $%
v_{i,j}\in \sigma _{i,j}(\overline{u}_{i})Z(H)^{m}\subset C_{H}(\overline{u}%
_{i}^{\ast })$.

For each $i\in \left\{ 1,\ldots ,n\right\} $, as the abelian group $C_{H}(%
\overline{u}_{i}^{\ast })/\Gamma (C_{H}(\overline{u}_{i}^{\ast }))$ is
finitely generated and torsion-free, it is freely generated by the images of
a sequence of elements $(w_{i,j})_{1\leq j\leq s(i)}$, where each $w_{i,j}$
can be written in the form $\tau _{i,j}(\overline{u}_{i})w^{m}$ with $w\in
C_{H}(\overline{u}_{i}^{\ast })$, which implies $w_{i,j}\in C_{H}(\overline{u%
}_{i}^{\ast })$.

For each $i\in \left\{ 1,\ldots ,n\right\} $, we denote by $H_{i}$ the
subgroup of $H$\ generated by $C_{H}(\overline{u}_{i}^{\ast })^{\prime }$, $%
\overline{\rho }_{i}(\overline{u}_{i})$, $v_{i,1},\ldots ,v_{i,r(i)}$, $%
w_{i,1},\ldots ,w_{i,s(i)}$. We have $H_{i}\subset C_{H}(\overline{u}%
_{i}^{\ast })$ and each element of $H_{i}$ can be written in a unique way as 
$uv_{i,1}^{a_{1}}\cdots v_{i,r(i)}^{a_{r(i)}}w_{i,1}^{b_{1}}\cdots
w_{i,s(i)}^{b_{s(i)}}$ with $u\in \overline{\rho }_{i}(\overline{u}%
_{i})C_{H}(\overline{u}_{i}^{\ast })^{\prime }$, $a_{1},\ldots ,a_{r(i)}\in 
%TCIMACRO{\U{2124} }%
%BeginExpansion
\mathbb{Z}
%EndExpansion
$ and $b_{1},\ldots ,b_{s(i)}\in 
%TCIMACRO{\U{2124} }%
%BeginExpansion
\mathbb{Z}
%EndExpansion
$.

We have $H=H_{1}\times \cdots \times H_{n}$ because of the following facts:

\noindent 1) For $i\neq j$, $[H_{i},H_{j}]=\left\{ 1\right\} $ follows from $%
[C_{H}(\overline{u}_{i}^{\ast }),C_{H}(\overline{u}_{j}^{\ast })]=\left\{
1\right\} $.

\noindent 2) We have $\Delta (H)=\overline{\rho }_{1}(\overline{u}_{1})C_{H}(%
\overline{u}_{1}^{\ast })^{\prime }\times \cdots \times \overline{\rho }_{n}(%
\overline{u}_{n})C_{H}(\overline{u}_{n}^{\ast })^{\prime }$.

\noindent 3) We have $\Gamma (H)/\Delta (H)=\langle v_{1,1}^{\prime },\ldots
,v_{1,r(1)}^{\prime }\rangle \times \cdots \times \langle v_{n,1}^{\prime
},\ldots ,v_{n,r(n)}^{\prime }\rangle $ where the $v_{i,j}^{\prime }$ are
the images of the $v_{i,j}$ in $H/\Delta (H)$.

\noindent 4) The properties $Z(H)=Z(C_{H}(\overline{u}_{1}^{\ast }))=\cdots
=Z(C_{H}(\overline{u}_{n}^{\ast }))$\ and

\noindent $H/Z(H)=C_{H}(\overline{u}_{1}^{\ast })/Z(H)\times \cdots \times
C_{H}(\overline{u}_{n}^{\ast })/Z(H)$ imply

\noindent $H/\Gamma (H)=C_{H}(\overline{u}_{1}^{\ast })/\Gamma (C_{H}(%
\overline{u}_{1}^{\ast }))\times \cdots \times C_{H}(\overline{u}_{n}^{\ast
})/\Gamma (C_{H}(\overline{u}_{n}^{\ast }))$

\noindent $=\langle w_{1,1}^{\prime },\ldots ,w_{1,s(1)}^{\prime }\rangle
\times \cdots \times \langle w_{n,1}^{\prime },\ldots ,w_{n,s(n)}^{\prime
}\rangle $

\noindent where the $w_{i,j}^{\prime }$ are the images of the $w_{i,j}$ in $%
C_{H}(\overline{u}_{i}^{\ast })/\Gamma (C_{H}(\overline{u}_{i}^{\ast }))$.

Now we prove that $\left| H_{i}:\mathrm{pr}_{H_{i}}(f(G_{i}))\right| $ is
prime to $m$ for $1\leq i\leq n$.

As $(\overline{u}_{1},\ldots ,\overline{u}_{n})$ satisfies $\gamma \wedge
\delta $ in $H$, $f$ induces some injective homomomorphisms from $G/\Gamma
(G)$ to $H/\Gamma (H)$ and from $\Gamma (G)/\Delta (G)$ to $\Gamma
(H)/\Delta (H)$. It follows $f(\Gamma (G))=\Gamma (H)\cap f(G)$ and $%
f(\Delta (G))=\Delta (H)\cap f(G)$.

Moreover, $\left| \Delta (H):f(\Delta (G))\right| $ is prime to $m$\ because
it divides $\left| H:f(G)\right| =\left| H:f(G)\Delta (H)\right| .\left|
\Delta (H):f(\Delta (G))\right| $. It follows that $\left| \Delta
(H_{i}):f(\Delta (G_{i}))\right| $ is prime to $m$ since $\Delta (G)=\Delta
(G_{1})\times \cdots \times \Delta (G_{n})$, $\Delta (H)=\Delta
(H_{1})\times \cdots \times \Delta (H_{n})$ and $f(\Delta (G_{j}))\subset
\Delta (H_{j})$ for $1\leq j\leq n$.

We have $v_{i,j}\mathrm{pr}_{H_{i}}(\sigma _{i,j}(\overline{u}_{i}))^{-1}=%
\mathrm{pr}_{H_{i}}(v_{i,j}\sigma _{i,j}(\overline{u}_{i})^{-1})\in \mathrm{%
pr}_{H_{i}}(H^{m})=H_{i}^{m}$ for $1\leq j\leq r(i)$ and $w_{i,j}\mathrm{pr}%
_{H_{i}}(\tau _{i,j}(\overline{u}_{i}))^{-1}=\mathrm{pr}_{H_{i}}(w_{i,j}\tau
_{i,j}(\overline{u}_{i})^{-1})\in \mathrm{pr}_{H_{i}}(H^{m})=H_{i}^{m}$ for $%
1\leq j\leq s(i)$. Consequently, $H_{i}=\left\langle \Delta
(H_{i}),(v_{i,j})_{1\leq j\leq r(i)},(w_{i,j})_{1\leq j\leq
s(i)}\right\rangle $ is generated by $\mathrm{pr}_{H_{i}}(f(G_{i}))$, $%
H_{i}^{m}$ and $\Delta (H_{i})$, and $\left| H_{i}:\mathrm{pr}%
_{H_{i}}(f(G_{i}))\Delta (H_{i})\right| $ is prime to $m$. It follows that $%
\left| H_{i}:\mathrm{pr}_{H_{i}}(f(G_{i}))\right| $ is prime to $m$ since $%
\mathrm{pr}_{H_{i}}(f(\Delta (G_{i})))=f(\Delta (G_{i}))$ and $\left| \Delta
(H_{i}):f(\Delta (G_{i}))\right| $ is prime to $m$.

Finally, we show that $H_{i}=\mathrm{pr}_{H_{i}}(f(G_{i}))Z(H_{i}^{r})$ for $%
1\leq i\leq n$.

We have $H^{\prime }\subset f(G)^{\prime }Z(H^{r})=f(G^{\prime })Z(H^{r})$
because $H=f(G)Z(H^{r})$\ and\ $Z(H^{r})$\ is normal in $H$. As $f(\Delta
(G))$\ contains each $\overline{\rho }_{j}(\overline{u}_{j})$, it follows $%
\Delta (H)\subset f(\Delta (G))Z(H^{r})$. Consequently, we have $\Delta
(H_{i})\subset f(\Delta (G_{i}))Z(H_{i}^{r})=\mathrm{pr}_{H_{i}}(f(\Delta
(G_{i})))Z(H_{i}^{r})$ since $\mathrm{pr}_{H_{i}}(\Delta (H))=\Delta (H_{i})$%
, $\mathrm{pr}_{H_{i}}(Z(H^{r}))=Z(H_{i}^{r})$, $\mathrm{pr}%
_{H_{i}}(f(\Delta (G_{i})))=f(\Delta (G_{i}))$\ and\ $\mathrm{pr}%
_{H_{i}}(f(\Delta (G_{j})))=\left\{ 1\right\} $ for $j\in \left\{ 1,\ldots
,n\right\} -\left\{ i\right\} $.

For $1\leq j\leq r(i)$, we have $v_{i,j}\in \sigma _{i,j}(\overline{u}%
_{i})Z(H)^{m}$, and therefore $v_{i,j}\in \mathrm{pr}_{H_{i}}(\sigma _{i,j}(%
\overline{u}_{i}))Z(H_{i})^{m}$ since $\mathrm{pr}_{H_{i}}(v_{i,j})=v_{i,j}$%
\ and $\mathrm{pr}_{H_{i}}(Z(H)^{m})=Z(H_{i})^{m}$. As $r$ divides $m$, we
have $Z(H_{i})^{m}\subset Z(H_{i})^{r}\subset Z(H_{i}^{r})$. Consequently,

\noindent $\mathrm{pr}_{H_{i}}(f(\Gamma (G_{i})))Z(H_{i}^{r})$\ contains $%
v_{i,1},\ldots ,v_{i,r(i)}$ and therefore contains\ $\Gamma (H_{i})$.

For each $z\in H_{i}$, as $H=f(G)Z(H^{r})$, there exist $x_{1}\in
\left\langle \overline{x}_{1}\right\rangle ,\ldots ,x_{n}\in \left\langle 
\overline{x}_{n}\right\rangle $, $y\in Z(H^{r})$ such that $z=f(x_{1})\cdots
f(x_{n})y$. It follows that $z=\mathrm{pr}_{H_{i}}(z)=\mathrm{pr}%
_{H_{i}}(f(x_{1}))\cdots \mathrm{pr}_{H_{i}}(f(x_{n}))\mathrm{pr}_{H_{i}}(y)$
belongs to $\mathrm{pr}_{H_{i}}(f(G_{i}))Z(H_{i}^{r})$\ since\ $\mathrm{pr}%
_{H_{i}}(y)\in Z(H_{i}^{r})$, $\mathrm{pr}_{H_{i}}(f(x_{i}))\in \mathrm{pr}%
_{H_{i}}(f(G_{i}))$\ and $\mathrm{pr}_{H_{i}}(f(x_{j}))\in Z(H_{i})\subset
\Gamma (H_{i})\subset \mathrm{pr}_{H_{i}}(f(\Gamma (G_{i})))Z(H_{i}^{r})$
for $j\in \left\{ 1,\ldots ,n\right\} -\left\{ i\right\} $.~~$\blacksquare $%
\bigskip

\noindent \textbf{Theorem 2.2.} Let $G,H$ be\ elementarily equivalent
polycyclic-by-finite groups. Then, for each decomposition $G\cong
G_{1}\times \cdots \times G_{m}$, there exists a decomposition $H\cong
H_{1}\times \cdots \times H_{m}$ with $H_{i}\equiv G_{i}$ for $1\leq i\leq m$%
.\bigskip

\noindent \textbf{Proof.} There exists $k\in 
%TCIMACRO{\U{2115} }%
%BeginExpansion
\mathbb{N}
%EndExpansion
^{\ast }$ such that $\Gamma (G)^{k}\subset Z(G)G^{\prime }$, $\Delta
(G)^{k}\subset G^{\prime }$, $\Gamma (H)^{k}\subset Z(H)H^{\prime }$ and $%
\Delta (H)^{k}\subset H^{\prime }$. Consequently, by Theorem 2.1, there
exist $r\in 
%TCIMACRO{\U{2115} }%
%BeginExpansion
\mathbb{N}
%EndExpansion
^{\ast }$ and, for each $n\in 
%TCIMACRO{\U{2115} }%
%BeginExpansion
\mathbb{N}
%EndExpansion
^{\ast }$, a decomposition $H\cong H_{1,n}\times \cdots \times H_{m,n}$ and
some injective homomorphisms $g_{i,n}:G_{i}\rightarrow H_{i,n}$ with $\left|
H_{i,n}/H_{i,n}^{r}\right| =\left| G_{i}/G_{i}^{r}\right| $, $%
H_{i,n}=g_{i,n}(G_{i})Z(H_{i,n}^{r})$ and $\left|
H_{i,n}:g_{i,n}(G_{i})\right| $\ prime to $n!$ for $1\leq i\leq m$.

According to [14, Cor. 3], up to isomorphism, $H$ only has finitely many
decompositions in direct products of indecomposable groups. It follows that,
up to isomorphism, $H$ only has finitely many decompositions in direct
products of groups, since any such decomposition is obtained by grouping
together some factors of a decomposition in direct product of indecomposable
groups. Consequently, there exists $S\subset 
%TCIMACRO{\U{2115} }%
%BeginExpansion
\mathbb{N}
%EndExpansion
^{\ast }$ infinite such that, for each $i\in \left\{ 1,\ldots ,m\right\} $,
the subgroups $H_{i,n}$ for $n\in S$ are all isomorphic. By Theorem 1.1, it
follows $G_{i}\equiv H_{i,n}$\ for $1\leq i\leq m$ and $n\in S$.~~$%
\blacksquare $\bigskip

\noindent \textbf{Remark.} In particular, if $H$ is indecomposable, then $G$
is indecomposable.\bigskip

The following definitions are slightly different from those given by F. Oger
in [14]. We are using them because they are easier to manage.\bigskip

\noindent \textbf{Definitions.} Denote by $J$ the infinite cyclic group with
multiplicative notation. Let $G$ be a group with $G/G^{\prime }$ finitely
generated. Then $G$ is $J$\emph{-indecomposable} if it is indecomposable and
if, for each $r\in 
%TCIMACRO{\U{2115} }%
%BeginExpansion
\mathbb{N}
%EndExpansion
^{\ast }$ and any groups $M,N$, $(\times ^{r}J)\times G\cong M\times N$
implies $M$ or $N$ torsion-free abelian (in particular, $J$ is $J$%
-indecomposable).\ A $J$\emph{-decomposition} of $G$ is a decomposition $%
(\times ^{r}J)\times G\cong G_{1}\times \cdots \times G_{n}$ with $%
G_{1},\ldots ,G_{n}$ $J$-indecomposable.\bigskip

According to [14, Prop. 1], each group $M$ with $M/M^{\prime }$ finitely
generated which satisfies the maximal condition on direct factors, and in
particular each polycyclic-by-finite group, has a $J$-decomposition. By [14,
Th.], for any $J$-indecomposable groups $G_{1},\ldots ,G_{m}$, $H_{1},\ldots
,H_{n}$ with $G_{1}/G_{1}^{\prime },\ldots ,G_{m}/G_{m}^{\prime }$, $%
H_{1}/H_{1}^{\prime },\ldots ,H_{n}/H_{n}^{\prime }$ finitely generated, if $%
G_{1}\times \cdots \times G_{m}\cong H_{1}\times \cdots \times H_{n}$, then $%
m=n$ and there exists a permutation $\sigma $ of $\left\{ 1,\ldots
,n\right\} $ such that $J\times G_{i}\cong J\times H_{\sigma (i)}$ for $%
1\leq i\leq n$.

According to [10], for any groups $G,H$, if $J\times G\cong J\times H$, then 
$G\equiv H$. By [14, Lemma 1], $J\times J\times G\cong J\times H$ implies $%
J\times G\cong H$.\bigskip

\noindent \textbf{Lemma 2.3.} For any elementarily equivalent
polycyclic-by-finite groups $G,H$, if $G$ is $J$-indecomposable, then $H$ is 
$J$-indecomposable.\bigskip

\noindent \textbf{Proof.} Otherwise, there exists a decomposition $(\times
^{r}J)\times H\cong H_{1}\times H_{2}$ with $H_{1},H_{2}$\ not torsion-free
abelian. As $(\times ^{r}J)\times G\equiv (\times ^{r}J)\times H$, Theorem
2.2 implies the existence of a decomposition $(\times ^{r}J)\times G\cong
G_{1}\times G_{2}$ with $G_{1}\equiv H_{1}$ and $G_{2}\equiv H_{2}$. It
follows that $G$ is not $J$-indecomposable.~~$\blacksquare $\bigskip

\noindent \textbf{Corollary 2.4.} For any\ elementarily equivalent
polycyclic-by-finite groups $G$, $H$, each $r\in 
%TCIMACRO{\U{2115} }%
%BeginExpansion
\mathbb{N}
%EndExpansion
^{\ast }$ and any $J$-decompositions $(\times ^{r}J)\times G\cong
G_{1}\times \cdots \times G_{m}$ and $(\times ^{r}J)\times H\cong
H_{1}\times \cdots \times H_{n}$, we have $m=n$ and there exists a
permutation $\sigma $ of $\left\{ 1,\ldots ,n\right\} $ such that $%
G_{i}\equiv H_{\sigma (i)}$ for $1\leq i\leq n$.\bigskip

\noindent \textbf{Proof.} As $(\times ^{r}J)\times G\equiv (\times
^{r}J)\times H$, Theorem 2.2 implies the existence of a decomposition $%
(\times ^{r}J)\times H\cong K_{1}\times \cdots \times K_{m}$ with $%
K_{i}\equiv G_{i}$ for $1\leq i\leq m$. By Lemma 2.3, each $K_{i}$\ is $J$%
-indecomposable. According to [14, Th.], we have $m=n$ and there exists a
permutation $\sigma $ of $\left\{ 1,\ldots ,n\right\} $ such that $J\times
H_{\sigma (i)}\cong J\times K_{i}$, and therefore $H_{\sigma (i)}\equiv
K_{i}\equiv G_{i}$, for $1\leq i\leq n$.~~$\blacksquare $\bigskip

\noindent \textbf{Lemma 2.5.} For any polycyclic-by-finite groups $G,H$ and
each $r\in 
%TCIMACRO{\U{2115} }%
%BeginExpansion
\mathbb{N}
%EndExpansion
^{\ast }$, if $(\times ^{r}J)\times G\equiv (\times ^{r}J)\times H$, then $%
G\equiv H$.\bigskip

\noindent \textbf{Proof.} It follows from Theorem 2.2 applied to $(\times
^{r}J)\times G$ and $(\times ^{r}J)\times H$ that there exists a group $%
K\equiv G$ such that $(\times ^{r}J)\times H\cong (\times ^{r}J)\times K$,
which implies $J\times H\cong J\times K$ and $H\equiv K\equiv G$.~~$%
\blacksquare $\bigskip

\noindent \textbf{Corollary 2.6.} For any polycyclic-by-finite groups $%
G_{1},G_{2},H_{1},H_{2}$ such that $G_{1}\equiv H_{1}$, we have $G_{1}\times
G_{2}\equiv H_{1}\times H_{2}$ if and only if $G_{2}\equiv H_{2}$.\bigskip

\noindent \textbf{Proof.} The condition is sufficient by Feferman-Vaught's
theorem. It remains to be proved that it is necessary.

By Corollary 2.4, there exist some $J$-decompositions $(\times ^{r}J)\times
G_{1}\cong G_{1,1}\times \cdots \times G_{1,k}$ and $(\times ^{r}J)\times
H_{1}\cong H_{1,1}\times \cdots \times H_{1,k}$ with $G_{1,i}\equiv H_{1,i}$
for $1\leq i\leq k$. We also consider some $J$-decompositions $(\times
^{s}J)\times G_{2}\cong G_{2,1}\times \cdots \times G_{2,m}$ and $(\times
^{s}J)\times H_{2}\cong H_{2,1}\times \cdots \times H_{2,n}$.

The property $G_{1}\times G_{2}\equiv H_{1}\times H_{2}$\ implies $(\times
^{r+s}J)\times G_{1}\times G_{2}\equiv (\times ^{r+s}J)\times H_{1}\times
H_{2}$. It follows from Corollary 2.4 applied to the $J$-decompositions $%
(\times ^{r+s}J)\times G_{1}\times G_{2}\cong G_{1,1}\times \cdots \times
G_{1,k}\times G_{2,1}\times \cdots \times G_{2,m}$ and $(\times
^{r+s}J)\times H_{1}\times H_{2}\cong H_{1,1}\times \cdots \times
H_{1,k}\times H_{2,1}\times \cdots \times H_{2,n}$ that, for each group $K$,
we have $\left| \left\{ i\in \left\{ 1,\ldots ,m\right\} \mid G_{2,i}\equiv
K\right\} \right| =\left| \left\{ j\in \left\{ 1,\ldots ,n\right\} \mid
H_{2,j}\equiv K\right\} \right| $ since $\left| \left\{ i\in \left\{
1,\ldots ,k\right\} \mid G_{1,i}\equiv K\right\} \right| =\left| \left\{
j\in \left\{ 1,\ldots ,k\right\} \mid H_{1,j}\equiv K\right\} \right| $.
Consequently, we have $(\times ^{s}J)\times G_{2}\equiv (\times ^{s}J)\times
H_{2}$. It follows $G_{2}\equiv H_{2}$\ by Lemma 2.5.~~$\blacksquare $%
\bigskip

\noindent \textbf{Corollary 2.7.} For any polycyclic-by-finite groups $G,H$
and each integer $k\geq 2$, we have $\times ^{k}G\equiv \times ^{k}H$ if and
only if $G\equiv H$.\bigskip

\noindent \textbf{Proof.} The condition is sufficient by Feferman-Vaught's
theorem. It remains to be proved that it is necessary.

We consider some $J$-decompositions $(\times ^{r}J)\times G\cong G_{1}\times
\cdots \times G_{m}$ and $(\times ^{r}J)\times H\cong H_{1}\times \cdots
\times H_{n}$. Then $(\times ^{kr}J)\times (\times ^{k}G)\cong (\times
^{k}G_{1})\times \cdots \times (\times ^{k}G_{m})$ is a $J$-decomposition of 
$\times ^{k}G$ and $(\times ^{kr}J)\times (\times ^{k}H)\cong (\times
^{k}H_{1})\times \cdots \times (\times ^{k}H_{n})$ is a $J$-decomposition of 
$\times ^{k}H$.

The property $\times ^{k}G\equiv \times ^{k}H$\ implies $(\times
^{kr}J)\times (\times ^{k}G)\equiv (\times ^{kr}J)\times (\times ^{k}H)$. By
Corollary 2.4, we have $m=n$ and, for each group $L$, the decompositions $%
(\times ^{k}G_{1})\times \cdots \times (\times ^{k}G_{n})$ and $(\times
^{k}H_{1})\times \cdots \times (\times ^{k}H_{n})$\ have the same number of
factors which are elementarily equivalent to $L$. Consequently, the same
property is true for the decompositions $G_{1}\times \cdots \times G_{n}$
and $H_{1}\times \cdots \times H_{n}$. It follows\ $G_{1}\times \cdots
\times G_{n}\equiv H_{1}\times \cdots \times H_{n}$ and $(\times
^{r}J)\times G\equiv (\times ^{r}J)\times H$. This property implies $G\equiv
H$\ by Lemma 2.5.~~$\blacksquare $\bigskip

\bigskip

\begin{center}
\textbf{References}\bigskip
\end{center}

\noindent \lbrack 1] P.M. Cohn, Algebra, Second Edition, Volume 2, John
Wiley \& Sons, 1989.

\noindent \lbrack 2] P.C. Eklof and E.R. Fisher, The elementary theory of
abelian groups, Ann. Math. Logic 4 (1972), 115-171.

\noindent \lbrack 3] S. Feferman and R.L. Vaught, The first order properties
of algebraic systems, Fund. Math. 47 (1959), 57-103.

\noindent \lbrack 4] F.J. Gr\"{u}newald, P.F. Pickel and D. Segal,
Polycyclic groups with isomorphic finite quotients, Annals of Math. 111
(1980), 155-195.

\noindent \lbrack 5] O. Kharlampovich and A. Myasnikov, Elementary theory of
free non-abelian groups, J. Algebra 302 (2006), 451-652.

\noindent \lbrack 6] C. Lasserre, Polycyclic-by-finite groups and
first-order sentences, J. Algebra 396 (2013), 18-38.

\noindent \lbrack 7] L.M. Manevitz, Applied model theory and
metamathematics. An Abraham Robinson memorial problem list, Israel J. Math.
49 (1984), 3-14.

\noindent \lbrack 8] A.G. Myasnikov and V.N. Remeslennikov, Model-theoretic
questions of group theory, Vopr. Algebry (Minsk Univ.) 4 (1989), 16-22
(Russian).

\noindent \lbrack 9] F. Oger, Equivalence \'{e}l\'{e}mentaire entre groupes
finis-par-abeliens de type fini, Comment Math. Helv. 57 (1982), 469-480
(French).

\noindent \lbrack 10] F. Oger, Cancellation and elementary equivalence of
groups, J. Pure Appl. Algebra 30 (1983), 293-299.

\noindent \lbrack 11] F. Oger, Elementary equivalence and profinite
completions: a characterization of finitely generated abelian-by-finite
groups, Proc. Amer. Math. Soc. 103 (1988), 1041-1048.

\noindent \lbrack 12] F. Oger, Cancellation and elementary equivalence of
finitely generated finite-by-nilpotent groups, J. London Math. Soc. (2) 44
(1991), 173-183.

\noindent \lbrack 13] F. Oger, Finite images and elementary equivalence of
completely regular inverse semigroups, Semigroup Forum 45 (1992), 322-331.

\noindent \lbrack 14] F. Oger, The direct decompositions of a group G with
G/G' finitely generated, Trans. Amer. Math. Soc 347 (1995), 1997-2010.

\noindent \lbrack 15] V.A. Roman'kov, Width of verbal subgroups in solvable
groups, Algebra and Logic 21 (1982), 41-49.

\noindent \lbrack 16] D. Segal, Words - Notes on Verbal Width in Groups,
London Math. Soc. Lecture Note Series 361, Cambridge Univ. Press, Cambridge,
2009.

\noindent \lbrack 17] Z. Sela, Diophantine geometry over groups VI, The
elementary theory of a free group, Geom. Funct. Anal. 16 (3) (2006),
707--730.

\noindent \lbrack 18] B.I. Zil'ber, An example of two elementarily
equivalent, but non isomorphic finitely generated metabelian groups, Algebra
I Logika 10 (1971), 309-315.\bigskip

\vfill

\end{document}